\documentclass[12pt,a4paper]{amsart}
\usepackage{amsmath,amstext,amssymb,amsthm,amsfonts}
\newcommand{\nc}{\newcommand}

\nc{\one}{\mbox{\bf 1}}
\nc{\invtensor}{\underset{\leftarrow}{\otimes}}
\nc{\const}{\operatorname{const}}

\nc{\ad}{\operatorname{ad}}
\nc{\tr}{\operatorname{tr}}
\nc{\tp}{\operatorname{top}}
\nc{\rank}{\operatorname{rank}}
\nc{\corank}{\operatorname{corank}}
\nc{\codim}{\operatorname{codim}}
\nc{\sdim}{\operatorname{sdim}}
\nc{\mult}{\operatorname{mult}}

\nc{\spn}{\operatorname{span}}
\nc{\Sym}{\operatorname{Sym}}
\nc{\sym}{\operatorname{sym}}
\nc{\id}{\operatorname{id}}
\nc{\Id}{\operatorname{Id}}
\nc{\Ree}{\operatorname{Re}}

\nc{\htt}{\operatorname{ht}}
\nc{\Ker}{\operatorname{Ker}}
\nc{\rker}{\operatorname{rKer}}
\nc{\im}{\operatorname{Im}}
\nc{\osp}{\mathfrak{osp}}
\nc{\sgn}{\operatorname{sgn}}
\nc{\F}{\operatorname{F}}
\nc{\Mod}{\operatorname{Mod}}
\nc{\Mat}{\operatorname{Mat}}
\nc{\Soc}{\operatorname{Soc}}
\nc{\Inj}{\operatorname{Inj}}
\nc{\Hom}{\operatorname{Hom}}
\nc{\End}{\operatorname{End}}
\nc{\supp}{\operatorname{supp}}
\nc{\Card}{\operatorname{Card}}
\nc{\Ann}{\operatorname{Ann}}
\nc{\Ind}{\operatorname{Ind}}
\nc{\Coind}{\operatorname{Coind}}
\nc{\wt}{\operatorname{wt}}
\nc{\ch}{\operatorname{ch}}
\nc{\Stab}{\operatorname{Stab}}
\nc{\Sch}{{\mathcal S}\mbox{\em ch}}
\nc{\Irr}{\operatorname{Irr}}
\nc{\Spec}{\operatorname{Spec}}
\nc{\Prim}{\operatorname{Prim}}
\nc{\Aut}{\operatorname{Aut}}
\nc{\Ext}{\operatorname{Ext}}

\nc{\Fract}{\operatorname{Fract}}
\nc{\gr}{\operatorname{gr}}
\nc{\deff}{\operatorname{def}}
\nc{\HC}{\operatorname{HC}}

\nc{\red}{\operatorname{red}}

\nc{\wdchi}{\widetilde{\chi}}
\nc{\wdH}{\widetilde{H}}
\nc{\wdN}{\widetilde{N}}
\nc{\wdM}{\widetilde{M}}
\nc{\wdO}{\widetilde{O}}
\nc{\wdR}{\widetilde{R}}
\nc{\wdS}{\widetilde{S}}
\nc{\wdV}{\widetilde{V}}

\nc{\wdC}{\widetilde{C}}

\nc{\Obj}{\operatorname{Obj}}
\nc{\Dglie}{\operatorname{{\mathcal D}glie}}
\nc{\Fin}{\operatorname{{\mathcal F}in}}

\nc{\Adm}{\operatorname{\mathcal{A}dm}}

\nc{\Sg}{{\cS(\fg)}}
\nc{\Shg}{{\cS(\fhg)}}
\nc{\Ug}{{\cU(\fg)}}
\nc{\Uhg}{{\cU(\fhg)}}
\nc{\Sh}{{\cS(\fh)}}
\nc{\Uh}{{\cU(\fh)}}
\nc{\Uhh}{{\cU(\fhh)}}
\nc{\Zg}{{{\mathcal{Z}}(\fg)}}

\nc{\Vir}{{\mathcal{V}ir}}

\nc{\tZg}{{\widetilde{\mathcal Z}({\mathfrak g})}}
\nc{\Zk}{{\mathcal Z}({\mathfrak k})}
\nc{\str}{\mathfrak{str}}

\nc{\Up}{{\mathcal U}({\mathfrak p})}
\nc{\Ah}{{\mathcal A}({\mathfrak h})}
\nc{\Ag}{{\mathcal A}({\mathfrak g})}
\nc{\Ap}{{\mathcal A}({\mathfrak p})}
\nc{\Zp}{{\mathcal Z}({\mathfrak p})}
\nc{\cZ}{\mathcal Z}
\nc{\cS}{\mathcal S}
\nc{\cT}{\mathcal{T}}
\nc{\cY}{\mathcal Y}
\nc{\cA}{\mathcal A}
\nc{\cU}{\mathcal U}
\nc{\cH}{\mathcal H}
\nc{\cM}{\mathcal M}
\nc{\cL}{\mathcal L}
\nc{\cF}{\mathcal F}
\nc{\fg}{\mathfrak g}

\nc{\fo}{\mathfrak o}

\nc{\CO}{\mathcal O}
\nc{\Cl}{\mathcal {C}\ell}

\nc{\cR}{\mathcal{R}}
\nc{\bM}{\mathbf{M}}
\nc{\bL}{\mathbf{L}}
\nc{\bN}{\mathbf{N}}

\nc{\zq}{\mathpzc q}

\nc{\fl}{\mathfrak l}
\nc{\fn}{\mathfrak n}
\nc{\fm}{\mathfrak m}
\nc{\fp}{\mathfrak p}
\nc{\fh}{\mathfrak h}
\nc{\ft}{\mathfrak t}
\nc{\fk}{\mathfrak k}
\nc{\fb}{\mathfrak b}
\nc{\fs}{\mathfrak s}
\nc{\fB}{\mathfrak B}

\nc{\vareps}{\varepsilon}
\nc{\varesp}{\varepsilon}
\nc{\veps}{\varepsilon}

\nc{\fsl}{\mathfrak{sl}}
\nc{\fpsl}{\mathfrak{psl}}
\nc{\fgl}{\mathfrak{gl}}
\nc{\fso}{\mathfrak{so}}

\nc{\fpq}{\mathfrak{pq}}
\nc{\fq}{\mathfrak q}
\nc{\fsq}{\mathfrak{sq}}
\nc{\fpsq}{\mathfrak{psq}}

\nc{\fhg}{\hat{\fg}}
\nc{\fhn}{\hat{\fn}}
\nc{\fhh}{\hat{\fh}}
\nc{\fhb}{\hat{\fb}}
\nc{\hrho}{\hat{\rho}}

\nc{\hsl}{\hat{\fsl}}
\nc{\fpo}{\mathfrak{po}}
\nc{\dirlim}{\underset{\rightarrow}{\lim}\,}
\nc{\nen}{\newenvironment}
\nc{\ol}{\overline}
\nc{\ul}{\underline}
\nc{\ra}{\rightarrow}
\nc{\lra}{\longrightarrow}
\nc{\Lra}{\Longrightarrow}

\nc{\Lla}{\Longleftarrow}

\nc{\Llra}{\Longleftrightarrow}

\nc{\thla}{\twoheadleftarrow}

\nc{\hra}{\hookrightarrow}

\nc{\iso}{\overset{\sim}{\lra}}

\nc{\ssubset}{\underset{\not=}{\subset}}

\nc{\vac}{|0\rangle}

\nc{\Thm}[1]{Theorem~\ref{#1}}
\nc{\Prop}[1]{Proposition~\ref{#1}}
\nc{\Lem}[1]{Lemma~\ref{#1}}
\nc{\Cor}[1]{Corollary~\ref{#1}}
\nc{\Conj}[1]{Conjecture~\ref{#1}}
\nc{\Claim}[1]{Claim~\ref{#1}}
\nc{\Defn}[1]{Definition~\ref{#1}}
\nc{\Exa}[1]{Example~\ref{#1}}
\nc{\Rem}[1]{Remark~\ref{#1}}
\nc{\Note}[1]{Note~\ref{#1}}
\nc{\Quest}[1]{Question~\ref{#1}}
\nc{\Hyp}[1]{Hypoth\`ese~\ref{#1}}
\nen{thm}[1]{\label{#1}{\bf Theorem.\ } \em}{}
\nen{prop}[1]{\label{#1}{\bf Proposition.\ } \em}{}
\nen{lem}[1]{\label{#1}{\bf Lemma.\ } \em}{}
\nen{cor}[1]{\label{#1}{\bf Corollary.\ } \em}{}
\nen{conj}[1]{\label{#1}{\bf Conjecture.\ } \em}{}

\nen{claim}[1]{\label{#1}{\bf Claim.\ } \em}{}

\nen{defn}[1]{\label{#1}{\bf Definition.\ } }{}
\nen{exa}[1]{\label{#1}{\bf Example.\ } }{}
\nen{rem}[1]{\label{#1}{\em Remark.\ } }{}
\nen{note}[1]{\label{#1}{\em Note.\ } }{}
\nen{exer}[1]{\label{#1}{\em Exercise.\ } }{}
\nen{sket}[1]{\label{#1}{\em Sketch of proof.\ } }{}
\nen{quest}[1]{\label{#1}{\bf Question.\ } \em}{}

\nen{hyp}[1]{\label{#1}{\bf Hypoth\`ese.\ } \em}{}
\begin{document}
\setcounter{section}{-1}

\title[Denominator identity]{Denominator identity
for affine Lie superalgebras with zero dual Coxeter number}
\author[Maria Gorelik, Shifra Reif]{Maria Gorelik, Shifra Reif}

\address{Dept. of Mathematics, The Weizmann Institute of Science,
Rehovot 76100, Israel}
\email{maria.gorelik@weizmann.ac.il, shifra.reif@weizmann.ac.i}
\thanks{Supported in part by ISF Grant No. 1142/07}

\begin{abstract}
We prove a denominator identity for non-twisted affine Lie superalgebras
with zero dual Coxeter number.
\end{abstract}

\maketitle
\section{Introduction}
\subsection{}
Let $\fg$ be a complex finite-dimensional contragredient Lie superalgebra. 
 These algebras were classified by V.~Kac in~\cite{Ksuper} 
and the list (excluding Lie algebras) consists of four series: 
$A(m|n), B(m|n), C(m), D(m|n)$ 
and the exceptional algebras $D(2,1,a), F(4), G(3)$.
The  finite-dimensional contragredient Lie superalgebras with 
zero Killing form (or, equivalently, with dual Coxeter number equal to zero) 
are $A(n|n), D(n|n+1)$ and $D(2,1,a)$.

Denote by $\Delta_{+0}$ (resp., $\Delta_{+1}$)
the set of positive even (resp., odd)  roots  of $\fg$.  
The Weyl denominator $R$ and the affine Weyl denominator 
$\hat{R}$ are given by the following formulas
$$R=\frac{R_0}{R_1},\ \ \hat{R}=\frac{\hat{R}_0}{\hat{R}_1},$$
where
$$\begin{array}{l}
R_0:=\prod_{\alpha\in\Delta_{+0}}(1-e^{-\alpha}),\ \ \
\hat{R}_0:=R_0\cdot\prod_{k=1}^{\infty}(1-q^k)^{\rank\fg}
\prod_{\alpha\in\Delta_{0}}(1-q^ke^{-\alpha}),\\
R_1:=\prod_{\alpha\in\Delta_{+1}}(1+e^{-\alpha}),\ \ \ 
\hat{R}_1:=R_1\cdot\prod_{k=1}^{\infty}
\prod_{\alpha\in\Delta_{1}}(1+q^ke^{-\alpha}).
\end{array}$$

Let $\hat{\fg}$ be the non-twisted affinization of $\fg$, $\fhh$
be the Cartan subalgebra of $\fhg$  and $\hat{\Delta}_+$
be the set of positive roots of $\fhg$.
The affine Weyl denominator is the Weyl denominator of $\fhg$.
Let $\hat{\rho}\in\fhh$ be such that  $2(\hat{\rho},\alpha)=(\alpha,\alpha)$
for each simple root $\alpha\in \hat{\Delta}_+$.

If $\fg$ has a non-zero Killing form, the affine denominator identity,
stated in~\cite{KW} and proven in~\cite{KW},\cite{G}, takes the form
\begin{equation}\label{denomnon0}
\hat{R}e^{\hat{\rho}}=\sum_{w\in T'} w(Re^{\hat{\rho}}),
\end{equation} 
where $T'$ is  the affine translation group corresponding
to the {\em ``largest''} root subsystem of $\Delta_0$ 
(see Section~\ref{Delta'non0} below). 
The affine denominator identity for strange Lie superalgebras $Q(n)$, 
which are not contragredient,
was stated in~\cite{KW} and proven in~\cite{Z}.

Suppose $\fg$ has zero dual Coxeter number, 
that is $\fg$ is $A(n|n),D(n|n+1)$ or $D(2,1,a)$.
In this case, 
$\hat{\rho}=\rho=\frac{1}{2}(\sum_{\alpha\in\Delta_{+0}}\alpha-
\sum_{\alpha\in\Delta_{+1}}\alpha)$.
In this paper we  will prove the following formulas
\begin{equation}\label{denom0}\begin{array}{ll}
\hat{R}e^{\hat{\rho}}\cdot f(q,e^{\str})=\sum_{w\in T'} w(Re^{\hat{\rho}})
& \text{ for }A(n|n),\\
\hat{R}e^{\hat{\rho}}\cdot f(q)=\sum_{w\in T'} w(Re^{\hat{\rho}})& 
\text{ for } D(n+1|n), D(2,1,a),
\end{array}\end{equation}
where $T'$ is the affine translation group corresponding
to the {\em ``smallest''} root subsystem of $\Delta_0$ (see~\ref{T'} below)
and $f(q,e^{\str}),f(q)$ are given by the formulas~(\ref{fqintr}) below.
The affine denominator identity
for $\fgl(2|2)$ was stated by V.~Kac and M.~Wakimoto in~\cite{KW}
and proven in~\cite{Ggl2} (the proof 
in~\cite{Ggl2} is different from the proof presented below).

In order to write down $f(q)$, we introduce the 
following infinite products after~\cite{DK}: 
for a parameter $q$ and a formal variable $x$ we set
$$(1+x)_q^{\infty}:=\prod_{k=0}^{\infty}(1+q^kx),\ \text{ and } \ 
(1-x)_q^{\infty}:=\prod_{k=0}^{\infty}(1-q^kx).$$
These infinite products converge
for any $x\in\mathbb{C}$ if the parameter $q$ is a real number $0<q<1$.
In particular, they are well defined for $0<x=q<1$ and 
 $(1\pm q)_q^{\infty}:=\prod_{n=1}^{\infty}(1\pm q^n)$. 

For $A(n|n)=\fgl(n|n)$ denote by  $\str$ the restriction of
the supertrace to the Cartan subalgebra $\fh\subset\fg$ 
(thus $\str\in\fh^*$). One has
\begin{equation}\label{fqintr}
\begin{array}{ll} 
f(q,e^{\str})=\frac{(1-q(-1)^ne^{\str})^{\infty}_{q}
\cdot(1-q(-1)^ne^{-\str})^{\infty}_{q}}{((1-q)^{\infty}_{q})^2}
& \text{ for } \fgl(n|n),\\
f(q)=\bigl((1-q)^{\infty}_{q}\bigr)^{-1}& \text{ for } D(n+1|n).
\end{array}
\end{equation}

As it was pointed by P.~Etingof, the terms $f(q,e^{\str}), f(q)$ 
can be interpreted using ``degenerate'' cases $n=1$; for example,
for $\fgl(1|1)$ we obtain the formula
$$\hat{R}e^{\hat{\rho}}=\frac{((1-q)^{\infty}_{q})^2}
{(1+qe^{\str})^{\infty}_{q}
\cdot(1+qe^{-\str})^{\infty}_{q}} Re^{\hat{\rho}},$$
which is trivial since $\fgl(1|1)$ has the only positive root $\beta=\str$,
which is odd.

Since $\fsl(n|n)=\{a\in\fgl(n|n)|\ \str(a) =0\}$ and $\rank \fsl(n|n)=2n-1=
\rank \fgl(n|n)-1$, one has  
$$f(q)=\left\{\begin{array}{ll} 
(1-q)^{\infty}_{q} & \text{ for } \fsl(2n|2n),\\
\frac{((1+q)^{\infty}_{q})^2}{(1-q)^{\infty}_{q}}
& \text{ for } \fsl(2n+1|2n+1).
\end{array}\right.$$

The root datum of $D(2,1,a)$ is the same as the root datum of $D(2|1)$
so the affine denominator identity for $D(2,1,a)$ is the same as 
the affine denominator identity for $D(2|1)$.

As it is shown in~\cite{KW}, the evaluation of
the affine denominator identity for $\fgl(2|2)$ 
(i.e., (\ref{denom0}) for $A(1|1)$)
gives the following Jacobi identity~\cite{J}:
\begin{equation}\label{jacobi}
\square(q)^8=1+16\sum_{j,k=1}^{\infty} (-1)^{(j+1)k}k^3q^{jk},
\end{equation}
where $\square(q)=\sum_{j\in\mathbb{Z}}q^{j^2}$ and thus
the coefficient of $q^m$ in the power series expansion of $\square(q)^8$
is the number  of representation of a given integer 
as a sum of $8$ squares (taking into the account the order of summands).

\subsection{}\label{T'}
In order to define $T'$ for $A(n|n), D(n+1|n)$ 
we present the set of even roots in the form
$\Delta_0=\Delta'\coprod\Delta''$, where
$$\begin{array}{ll}
\Delta'\cong \Delta''=A_{n-1} &\text{ for } A(n-1|n-1)=\fgl(n|n),\\
\Delta'=C_n,\ \Delta''=D_{n+1} & \text{ for } D(n+1|n).
\end{array}$$
Let $W'$ be the Weyl group of $\Delta'$ and $\hat{W}'$  
be  the corresponding affine Weyl group. 
Then $\hat{W}'=W'\ltimes T'$, where $T'$ is a translation group,
see~\cite{Kbook}, Chapter VI. Notice that for $D(n+1|n)$ the rank of
root system $\Delta'$ is smaller than the rank of $\Delta''$;
by contrast, for 
Lie superalgebras with non-zero Killing form, the lattice
$T'$ in~(\ref{denomnon0}) corresponds to the root system
$\Delta'$, whose rank is not smaller than the rank of $\Delta''$
(one has $\Delta_0=\Delta'\coprod\Delta''$ as before). It is not possible
to change $T'$ to $T''$ in Identity~(\ref{denomnon0}) and
in Identity~(\ref{denom0}) for $D(n+1|n)$, since
the sum $\sum_{w\in T''} w(Re^{\rho})$ is not well defined 
if $\Delta'\not\cong\Delta''$ (see~\Rem{T''}).

We prove Identity~(\ref{denom0}) and outline a similar proof for
Identity~(\ref{denomnon0}).
The key point is~\Prop{ai}, where it is shown that 
for any complex finite-dimensional  contragredient Lie superalgebra,
the expansion of
$Y:=\hat{R}^{-1}e^{-\hat{\rho}}\sum_{w\in T'} w(Re^{\hat{\rho}})$ contains only 
$\hat{W}$-invariant elements. This implies that 
$Y=f(q)$ for $\fg\not=\fgl(n|n)$ and $Y=f(q,e^{-\str})$ for $\fgl(n|n)$. 
We determine $f(q)$ for $D(n+1|n)$ and
$f(q,e^{\str})$ for $\fgl(n|n)$ using suitable evaluations.
For other finite-dimensional  contragredient simple Lie superalgebras
the equality $f(q)=1$ can be obtained in two steps:
first, using the Casimir operator and the fact that the dual Coxeter
number is non-zero, we show that $f(q)$ is scalar;
then  one deduces that this scalar is equal to $1$ from the denominator
identity for $\fg$ (this is done in~\cite{G}).

{\em Acknowledgement.} A part of this manuscript was written during
 the first author's stay at Max Planck Institut f\"ur Mathematik in Bonn, 
whose hospitality is greatly appreciated. The authors are grateful
to P.~Etingof and to V.~Kac for fruitful discussions.

\section{Preliminary}\label{sect1}
One readily sees (for instance,~\cite{G}, 1.5) that 
$Re^{\hat{\rho}}$ and $\hat{R}e^{\hat{\rho}}$
do not depend on the choice of set of positive roots $\Delta_{+}$
so it is enough to establish the identity for one choice of
$\Delta_{+}$. Similarly, it is enough
to establish the identity for one choice of  $A_{n-1}$ for $\fgl(n|n)$.
In Section~\ref{rootsys} we describe our choice of the
set of of positive roots for $\fgl(n|n), D(n+1|n)$.
In Section~\ref{ntaff} we introduce notation for affine
Lie superalgebra $\fhg$. In Section~\ref{cR} we introduce
the algebra $\cR$ of formal power series in which we expand $R$
and $\hat{R}$.

\subsection{Root systems}\label{rootsys}
Let $\fg$ be $\fgl(n|n)$ or $D(n|n+1)$ and let  $\fh$ be 
its Cartan subalgebra. 
We fix the following  sets of simple roots:
$$\begin{array}{l}
\Pi=\{\vareps_1-\delta_1,\delta_1-\vareps_2,\vareps_2-\delta_2,\ldots,
\vareps_n-\delta_n\} \text{ for } \fgl(n|n),\\
\Pi=\{\vareps_1-\delta_1,\delta_1-\vareps_2,\vareps_2-\delta_2,\ldots,
\vareps_n-\delta_n,\delta_n\pm\vareps_{n+1}\} \text{ for } D(n+1|n).
\end{array}$$

We fix a non-degenerate symmetric invariant bilinear 
form on $\fg$ and denote by $(-,-)$ the induced non-degenerate symmetric
bilinear form on $\fh^*$; we normalize the form in such a way that
$-(\vareps_i,\vareps_j)=(\delta_i,\delta_j)=\delta_{ij}$; notice that
$\{\vareps_i,\delta_i| 1\leq i\leq n\}$ (resp., 
$\{\vareps_j,\delta_i| 1\leq i\leq n, 1\leq j\leq n+1\}$
is an orthogonal basis of $\fh^*$ for $\fgl(n|n)$ (resp., for $D(n+1|n)$).

For this choice one has
$$\begin{array}{l}
\Delta_{0+}=\{\vareps_i-\varesp_j\}_{1\leq i<j\leq n}\coprod
\{\delta_i-\delta_j\}_{1\leq i<j\leq n}\ \text{ for } \fgl(n|n),\\
\Delta_{1+}=\{\vareps_i-\delta_j\}_{1\leq i\leq j\leq n}\cup
\{\delta_i-\vareps_j\}_{1\leq i<j\leq n}\ 
\text{ for } \fgl(n|n),\\
\Delta_{0+}=\{\vareps_i\pm\varesp_j\}_{1\leq i<j\leq n+1}
\coprod
\{\delta_s\pm\delta_t\}_{1\leq s<t\leq n}\cup
\{2\delta_s\}_{1\leq s\leq n}\  
\text{ for } D(n+1|n),\\
\Delta_{1+}=\{\vareps_i-\delta_s\}_{1\leq i\leq s\leq n}\cup
\{\delta_s-\vareps_j\}_{1\leq s<j\leq n+1}\cup\{\delta_i+\vareps_j\}_{1\leq
  i\leq n; 1\leq j\leq n+1}\text{ for } D(n+1|n).
\end{array}$$

For  $D(n+1|n)$ one has $\rho=0$ for $D(n+1|n)$. For $\fgl(n|n)$ one has
$\str=\sum_{i=1}^n (\vareps_i-\delta_i)$ and $\rho=-\frac{1}{2}\str$.

Recall that $\fsl(n|n)=\{a\in\fgl(n|n)|\ \str(a)=0\}$ and so
$\fh^*$ for $\fsl(n|n)$ is the quotient of $\fh^*$ for $\fgl(n|n)$
by $\mathbb{C}\str$.

By above, $\Delta_0$ is the union of two irreducible root systems, and
we write $\Delta_0=\Delta''\coprod\Delta'$, where $\Delta''$ lies 
in the span of $\vareps_i$s and $\Delta'$ lies in the span of $\delta_i$s 
(this notation is compatible with notations in Section~\ref{T'}).

\subsection{Non-twisted affinization}\label{ntaff}
Let $\fg=\fn_-\oplus\fh\oplus\fn_+$ be any complex 
finite-dimensional contragredient
Lie superalgebra with a fixed triangular decomposition,
and let $\Delta_+$ be its set of positive roots. 
Let $\fhg$ be the affinization of $\fg$ and 
let $\fhh$ be its Cartan subalgebra, see~\cite{Kbook}, Chapter VI.
Recall that $\fg=[\fg,\fg]\oplus\mathbb{C}D$ for some $D\in\fhh$.
Let $\hat{\Delta}=\hat{\Delta}_0\coprod
\hat{\Delta}_1$ be the  set of roots of $\fhg$. We  set
$$\hat{\Delta}^+=\Delta_+\cup (\displaystyle\cup_{k=1}^{\infty} \{\alpha+
k\delta|\ \alpha\in\Delta\})\cup (\displaystyle\cup_{k=1}^{\infty}
\{k\delta\}),
$$
where $\delta$ is the minimal imaginary root.
Let $W$ (resp., $\hat{W}$)  be the  Weyl group
of $\Delta_0$ (resp., $\hat{\Delta}_0$).
One has  $(\fhh^*)^{\hat{W}}=\mathbb{C}\delta$ for $\fg\not=\fgl(n|n)$
and $(\fhh^*)^{\hat{W}}=\mathbb{C}\delta\oplus \mathbb{C}\str$
for $\fg=\fgl(n|n)$.

We extend the non-degenerate symmetric invariant
bilinear form from $\fg$ to $\fhg$ and denote by $(-,-)$ the induced
non-degenerate symmetric bilinear form on $\fhh^*$
(the above-mentioned form on $\fh^*$ is induced by this  form on $\fhh^*$).
For $A\subset \fhh^*$ we set $A^{\perp}=\{\mu\in\fhh^*|\ \forall\nu\in A
\ (\mu,\nu)=0\}$. 

\subsubsection{}\label{Delta'non0}
In Section~\ref{rootsys} we introduced the root systems
$\Delta',\Delta''$ for $\fg=\fgl(n|n), D(n+1|n)$. 
For $\fg\not=\fgl(n|n), D(n+1|n), D(2,1,a)$
the Killing form $\kappa$ is non-zero; in this case, we introduce
$\Delta',\Delta''$ by the formulas:
$\Delta':=\{\alpha| \kappa(\alpha,\alpha)>0\}$,
$\Delta'':=\{\alpha| \kappa(\alpha,\alpha)<0\}$. One has
$\Delta_0=\Delta'\coprod\Delta''$ and
 $\Delta''=\emptyset$ if $\Delta_0$ is irreducible. 
Let $W'$ (resp., $W''$) be the Weyl group of $\Delta'$ (resp., $\Delta''$).
One has $W=W'\times W''$.

\subsubsection{}\label{TT}
Now that we have introduced the decomposition $\Delta_0=\Delta'\coprod\Delta''$
for any complex finite-dimensional contragredient Lie superalgebra,
we denote by $\hat{W}'$ the  Weyl group of the affine root system
$\hat{\Delta}'$. Recall that $\hat{W}'=W'\ltimes T'$,
where $T'$ is a translation group (see~\cite{Kbook}, Chapter VI).

\subsubsection{}
For $N\subset \fhh^*$ we use the notation $\mathbb{Z}N$ for
the set $\sum_{\mu\in N}\mathbb{Z}\mu$. 
Set
$$Q^+:=\sum_{\mu\in \Delta_+}\mathbb{Z}_{\geq  0}\mu,\ \ 
Q:=\mathbb{Z}\Delta,\ \ \
\hat{Q}^{\pm}:=\pm\sum_{\mu\in\hat{\Delta}_+}\mathbb{Z}_{\geq 0}\mu,\ \ 
\hat{Q}:=\mathbb{Z}\hat{\Delta}_+.$$
We introduce the standard partial order on $\fhh^*$:
$\mu\leq\nu$ if $(\nu-\mu)\in \hat{Q}^+$.

\subsection{Algebra $\cR$}\label{cR}
We are going to use notation of~\cite{G}, 1.4, which we recall below.
Retain notation of Section~\ref{ntaff}.

\subsubsection{}\label{R6}
Call a {\em $\hat{Q}^+$-cone} a set of the form $(\lambda-\hat{Q}^+)$, 
where $\lambda\in\hat{\fh}^*$.

For a formal sum of the form $Y:=\sum_{\nu\in\hat{\fh}^*} b_{\nu} e^{\nu},\ 
b_{\nu}\in\mathbb{Q}$ define the {\em support} of $Y$ by
$\supp(Y):=\{\nu\in\hat{\fh}^*|\ b_{\nu}\not=0\}$.
Let $\cR$ be a vector space over $\mathbb{Q}$,
spanned by the sums of the form
$\sum_{\nu\in \hat{Q}^+} b_{\nu} e^{\lambda-\nu}$, 
where $\lambda\in\hat{\fh}^*,\ 
b_{\nu}\in\mathbb{Q}$. In other words, $\cR$ consists of
the formal sums $Y=\sum_{\nu\in\hat{\fh}^*} b_{\nu}e^{\nu}$ with the support 
lying in a finite union of $\hat{Q}^+$-cones.

Clearly, $\cR$ has a structure of commutative algebra over 
$\mathbb{Q}$. If $Y\in \cR$ is such that $YY'=1$ for some $Y'\in\cR$,
we write $Y^{-1}:=Y'$.

\subsubsection{Action of the Weyl group}
For $w\in \hat{W}$ set $w(\sum_{\nu\in\hat{\fh}^*} b_{\nu}e^{\nu}):=
\sum_{\nu\in\hat{\fh}^*} b_{\nu}e^{w\nu}$. By above, $wY\in\cR$ iff
$w(\supp Y)$ is a subset of a finite union of $\hat{Q}^+$-cones.
For each subgroup $\tilde{W}$ of $\hat{W}$
we set $\cR_{\tilde{W}}:=\{Y\in\cR|\ wY\in \cR \text{ for
each }w\in \tilde{W}\}$; notice that $\cR_{\tilde{W}}$ is a subalgebra of $\cR$.
 
\subsubsection{Infinite products}\label{infprod}
An infinite product of the form $Y=\prod_{\nu\in X}
(1+a_{\nu}e^{-\nu})^{r(\nu)}$, where $a_{\nu}\in \mathbb{Q},\
\ r(\nu)\in\mathbb{Z}_{\geq 0}$ and $X\subset \hat{\Delta}$ is such that 
the set $X\setminus\hat{\Delta}_+$ is finite, can be naturally viewed 
as an element of $\cR$; clearly, this element does not depend
on the order of factors. Let $\cY$ be the set of such infinite products.
For any $w\in \hat{W}$ the infinite product
$$wY:=\prod_{\nu\in X}(1+a_{\nu}e^{-w\nu})^{r(\nu)},$$
is again an infinite product of the above form, since 
the set $w\hat{\Delta}_+\setminus \hat{\Delta}_+$  is finite
(see for example~\cite{G}, Lemma 1.2.8). 
Hence $\cY$ is a $\hat{W}$-invariant multiplicative subset of
$\cR_{\hat{W}}$.
 
The elements of $\cY$ are invertible in $\cR$: using 
the geometric series we can expand  $Y^{-1}$ 
(for example, $(1-e^{\alpha})^{-1}=
-e^{-\alpha}(1-e^{-\alpha})^{-1}=-\sum_{i=1}^{\infty} e^{-i\alpha}$).

\subsubsection{The subalgebra $\cR'$}\label{cR'}
Denote by $\cR'$ the localization of $\cR_{\hat{W}}$ by $\cY$. By above,
$\cR'$ is a subalgebra of $\cR$. Observe that $\cR'\not\subset \cR_{\hat{W}}$:
for example, $(1-e^{-\alpha})^{-1}\in\cR'$, but 
$(1-e^{-\alpha})^{-1}=\sum_{j=0}^{\infty} e^{-j\alpha}\not\in \cR_{\hat{W}}$.
We extend the action of $\hat{W}$
from $\cR_{\hat{W}}$ to $\cR'$ by setting $w(Y^{-1}Y'):=(wY)^{-1}(wY')$
for $Y\in\cY,\ Y'\in\cR_{\hat{W}}$.

Notice that an infinite product of the form $Y=\prod_{\nu\in X}
(1+a_{\nu}e^{-\nu})^{r(\nu)}$, where $a_{\nu}, X$ are as above
and $r(\nu)\in\mathbb{Z}$, lies in $\cR'$ and
$wY=\prod_{\nu\in X} (1+a_{\nu}e^{-w\nu})^{r(\nu)}$.
The support $\supp(Y)$ has a unique maximal element (with respect to the
standard partial order) and this element is given by the formula
$$\max\supp(Y)=-\sum_{\nu\in X\setminus\hat{\Delta}_+: a_{\nu}\not=0} 
r_{\nu}\nu.$$

\subsubsection{}\label{compex}
Let $\tilde{W}$ be a subgroup of $\hat{W}$.
For $Y\in\cR'$ we say that {\em $Y$ is $\tilde{W}$-invariant
(resp., $\tilde{W}$-anti-invariant)} if $wY=Y$
(resp., $wY=\sgn(w)Y$) for each $w\in \tilde{W}$.

Let $Y=\sum a_{\mu} e^{\mu}\in\cR_{\tilde{W}}$ be $\tilde{W}$-anti-invariant. 
Then $a_{w\mu}=(-1)^{\sgn(w)}a_{\mu}$
for each $\mu$ and $w\in \tilde{W}$. In particular,  
$\tilde{W}\supp(Y)=\supp(Y)$, and, moreover, for each  
$\mu\in\supp(Y)$ one has 
$\Stab_{\tilde{W}}\mu\subset\{w\in \tilde{W}|\ \sgn(w)=1\}$.
The condition $Y\in \cR_{\tilde{W}}$ is essential: for example, for
$\tilde{W}=\{\id,s_{\alpha}\}$, the expressions $Y:=e^{\alpha}-e^{-\alpha}$,
$Y^{-1}=e^{-\alpha}(1-e^{-2\alpha})^{-1}$ are  $\tilde{W}$-anti-invariant,
$\supp(Y)=\{\pm\alpha\}$ is $s_{\alpha}$-invariant,
but $\supp(Y^{-1})=\{-\alpha,-3\alpha,\ldots\}$ is not 
$s_{\alpha}$-invariant.

For $Y\in\cR_{\tilde{W}}$ such that
each $\tilde{W}$-orbit in $\fhh^*$ has a finite intersection with $\supp(Y)$,
introduce the sum
$$\cF_{\tilde{W}}(Y):=\sum_{w\in \tilde{W}}\!\sgn(w) wY.$$
This sum is well defined, but does not
always belong to $\cR$. For $Y=\sum a_{\mu} e^{\mu}$ one has
$\cF_{\tilde{W}}(Y)=\sum b_{\mu}e^{\mu}$, where
$b_{\mu}=\sum_{w\in \tilde{W}}\sgn(w) a_{w\mu}$; in particular,
$b_{\mu}=\sgn(w)b_{w\mu}$ for each $w\in \tilde{W}$. One has
$$Y\in\cR_{\tilde{W}}\ \&\ 
\cF_{\tilde{W}}(Y)\in\cR\ \Longrightarrow\ \left\{\begin{array}{l}
\supp (\cF_{\tilde{W}}(Y))\text{ is  $\tilde{W}$-stable},\\
\cF_{\tilde{W}}(Y)\in \cR_{\tilde{W}};\\  \cF_{\tilde{W}}(Y)
\text{ is $\tilde{W}$-anti-invariant}.
\end{array}\right.$$

We call a vector $\lambda\in\hat{\fh}^*$ {\em $\tilde{W}$-regular} 
if $\Stab_{\tilde{W}} \lambda=\{\id\}$, and we say
that the orbit $\tilde{W}\lambda$ is  $\tilde{W}$-regular if 
$\lambda$ is $\tilde{W}$-regular
(so the orbit consists of  $\tilde{W}$-regular points). 
If $\tilde{W}$ is an affine Weyl group, then for any $\lambda\in\hat{\fh}^*$
the stabilizer $\Stab_{\tilde{W}} \lambda$ is either trivial or 
contains a reflection.
Thus for $\tilde{W}=\hat{W}',\ \hat{W}''$ one has 
$$Y\in\cR_{\tilde{W}}\ \&\ 
\cF_{\tilde{W}}(Y)\in\cR\ \Longrightarrow\ 
\supp (\cF_{\tilde{W}}(Y))\ \text{ is a union
of $\tilde{W}$-regular orbits}.$$

For $Y\in\cR'$ the sum $\sum_{w\in \tilde{W}}\!\sgn(w) wY$ is not always 
$\tilde{W}$-anti-invariant: for example, for $\tilde{W}=\{\id,s_{\alpha}\}$
 one has $\sum_{w\in \tilde{W}}\!\sgn(w) w((1-e^{-\alpha})^{-1})=
(1-e^{-\alpha})^{-1}-(1-e^{\alpha})^{-1}=1+2e^{-\alpha}+2e^{-2\alpha}+\ldots$,
which is not $\tilde{W}$-anti-invariant.

\section{Proof}
As it is pointed out in Section~\ref{sect1}, it is enough to establish the
denominator identity for a particular choice of $\Delta_+$ and we do this
for the choice  described in Section~\ref{rootsys}. Recall that the group
$T'$ was introduced in Section~\ref{TT}.
The steps of the proof are the following. 

1) In Section~\ref{welldefin} we check  that for $\fg=\fgl(n|n), D(n+1|n)$,
the sum $\cF_{T'}(Re^{\hat{\rho}})$ is well-defined and belongs to $\cR$. 

2) In Section~\ref{pfsuppU} we prove the inclusions 
\begin{equation}\label{suppUU}
\supp(\cF_{T'}(Re^{\hat{\rho}})),
\supp(\hat{R}e^{\hat{\rho}})\subset U,
\end{equation}
where
\begin{equation}\label{defU}
U:=\{\mu\in\hat{\rho}-\hat{Q}^+|\ (\mu,\mu)=
(\hat{\rho},\hat{\rho})\}
\end{equation}
for $\fg=\fgl(n|n)$ and  $D(n+1|n)$. 

For simple contragredient Lie superalgebras with non-zero Killing form
steps (1), (2) are performed in~\cite{G}, 2.4.

3) In Section~\ref{Y1Y2} we show that for any finite-dimensional 
simple contragredient Lie superalgebra $\fg$
the inclusions~(\ref{suppUU}) imply that
$\supp\bigl(\hat{R}^{-1}e^{-\hat{\rho}}\cF_{T'}(Re^{\hat{\rho}})\bigr)
\subset \hat{Q}^{\hat{W}}$. As a result, 
$\hat{R}^{-1}e^{-\hat{\rho}}\cF_{T'}(Re^{\hat{\rho}})$ takes
the form $f(q)$ (resp., $f(q,e^{\str})$)
for $\fg\not=\fgl(n|n)$ (resp., for $\fgl(n|n)$).

4) In Section~\ref{ev} we compute $f(q)$ 
(resp., $f(q,e^{\str})$) for $D(n+1|n)$
(resp., for $\fgl(n|n)$). This completes the proof of 
Identity~(\ref{denom0}).

In Section~\ref{other} we briefly repeat the arguments
of~\cite{G} showing that $f(q)=1$ for $\fg\not=\fgl(n|n), D(n+1|n), D(2,1,a)$.
This completes the proof of Identity~(\ref{denomnon0}).

\subsection{Step 1}\label{welldefin}
In this subsection we show that for $\fg=\fgl(n|n), D(n+1|n)$,
the sum $\cF_{T'}(Re^{\hat{\rho}})$ 
is a well-defined element of $\cR$. Since $\hat{\rho}=\rho$
is $\hat{W}$-invariant, it is enough to verify that
$\cF_{T'}(R)$ is a well-defined element of $\cR$.

Recall that $T'=\mathbb{Z}\{t_{\delta_i-\delta_{i+1}}\}_{i=1}^{n-1}$ 
for $\fgl(n|n)$ and $T'=\mathbb{Z}\{t_{\delta_i}\}_{i=1}^n$ for $D(n+1|n)$,
where 
\begin{equation}\label{tmu}
t_{\mu}(\alpha)=\alpha-(\alpha,\mu)\delta\ \text{ for any }\alpha\in\hat{Q}.
\end{equation}

\subsubsection{}\label{welldef}
By Section~\ref{cR'} one has
$$\max\supp\bigl(w(R)\bigr)=\sum_{\alpha\in\Delta_{0+}: w\beta<0}w\alpha
-\sum_{\beta\in \Delta_{1+}: w\beta<0}w\beta.$$

For $w\in T'$ write  $w=t_{\mu}$, where $\mu\in 
\mathbb{Z}\{\delta_i-\delta_{i+1}\}_{1\leq i<n}$ 
for $\fgl(n|n)$ and $\mu\in\mathbb{Z}\{\delta_i\}_{i=1}^n$ for $D(n+1|n)$.
From~(\ref{tmu}) we get
$$\{\beta\in\Delta_{i+}| w\beta<0\}=
\{\beta\in\Delta_{i+}| (\beta,\mu)>0\} \ \text{ for }i=0,1.$$
We obtain
$\ \ \max\supp\bigl(t_{\mu}(R)\bigr)=-v(\mu)+(v(\mu),\mu)\delta$, where
$$v(\mu):=\sum_{\beta\in\Delta_{0+}: (\beta,\mu)>0}\beta-
\sum_{\beta\in\Delta_{1+}: (\beta,\mu)>0}\beta.
$$

In order to prove that $\cF_{T'}(R)$ is a well-defined element of $\cR$
we  verify that
\begin{equation}\label{eqwelldef}
(i)\ \ \forall\mu\ (v(\mu),\mu)\leq 0;\ \ (ii)\ \ 
\forall N>0\ \ \{\mu|\ (v(\mu),\mu)\geq -N\}\ \text{ is finite}.
\end{equation}
The condition (ii) ensures that the sum $\cF_{T'}(R)=\sum_{\mu} t_{\mu}(R)$
is well-defined and the condition (i) means that for each $\mu$ one has
$$\max\supp(t_{\mu}(R))=-v(\mu)\leq \sum_{\beta\in\Delta_{1+}} \beta$$
so $\supp\bigl(\cF_{T'}(R)\bigr)\subset \sum_{\beta\in\Delta_{1+}} \beta
-\hat{Q}^+$ and thus  $\cF_{T'}(R)\in\cR$.

\subsubsection{Case $\fgl(n|n)$}
Recall that $w\in T'$ has the form  $w=t_{\mu},\ \mu=\sum_{i=1}^n k_i\delta_i$,
where the $k_i$s are integers and $\sum_{i=1}^nk_i=0$. One has
$$\begin{array}{l}
\{\alpha\in\Delta_{+0}|\ (\alpha,\mu)>0\}:=
\{\delta_i-\delta_j|\ i<j, k_i>k_j\},\\
\{\alpha\in\Delta_{+1}|\ (\alpha,\mu)>0\}:=\{\vareps_i-\delta_j|\ k_j<0,
i\leq j\}\cup\{\delta_i-\vareps_j|\ k_i>0, i<j\},
\end{array} $$
where $1\leq i,j\leq n$.

Write $v(\mu)=v'+v''$, where $v'=\sum_{i=1}^n a_i\delta_i$  and
$v''$ lies in the span of $\vareps_i$s. By above,
for $k_i>0$ one has $a_i\leq (n-i)-(n-i)=0$ and
for $k_j<0$ one has $a_j\geq  -(j-1)+j=1$.
Therefore $(v(\mu),\mu)=\sum_{i=1}^n a_ik_i\leq \sum_{k_i<0} k_i\leq 0$ 
and the set $\{\mu|\ (v(\mu),\mu)\geq -N\}$ is a subset of the set 
$\{\mu|\ \sum_{k_i<0} k_i\geq -N\}$, 
which is finite for any $N$, because 
$k_i$s are integers and $\sum_{i=1}^n k_i=0$.
This establishes conditions~(\ref{eqwelldef}).

\subsubsection{Case $D(n+1|n)$}\label{welldefDnn+1}
Recall that $w\in T'$ has the form  $w=t_{\mu},\ \mu=\sum k_i\delta_i$,
where the $k_i$s are  integers. One has
$$\begin{array}{l}
\{\alpha\in\Delta_{+0}|\ (\alpha,\mu)>0\}:=
\{\delta_i-\delta_j|\ i<j, k_i>k_j\}\cup 
\{\delta_i+\delta_j|\ i\not=j,\ k_i+k_j>0\}\cup\{2\delta_i|\ k_i>0\},\\
\{\alpha\in\Delta_{+1}|\ (\alpha,\mu)>0\}:=\{\vareps_s-\delta_j|\ k_j<0,
s\leq j\}\cup\{\delta_i-\vareps_s|\ k_i>0, i<s\}\cup
\{\delta_i+\vareps_s|\ k_i>0\},
\end{array} $$
where $1\leq i,j\leq n$ and $1\leq s\leq n+1$.

Write $v(\mu)=v'+v''$, where $v'=\sum_{i=1}^n a_i\delta_i$  and
$v''$ lies in the span of $\vareps_i$s. By above,
for $k_i>0$ one has $a_i\leq (2n+1-i)-(2n+2-i)=-1$ and
for $k_j<0$ one has $a_j\geq  -(j-1)+j=1$.
Therefore 
$$(v(\mu),\mu)=\sum_{i=1}^n a_ik_i
\leq -\sum_{k_i>0} k_i+\sum_{k_j<0} k_j=-\sum_{1=1}^n |k_i|\leq 0$$ 
so the set
$\{\mu|\ (v(\mu),\mu)\geq -N\}$ is a subset of the set 
$\{\mu|\ \sum_{i=1}^n |k_i|\leq N\}$, which is finite for any $N$.
This establishes conditions~(\ref{eqwelldef}).

\subsubsection{}
\begin{rem}{T''}
For $\fgl(n|n)$ one can interchange $\Delta'$ and $\Delta''$ so
the sum  $\cF_{T''}(R)$ is well-defined. 
One readily sees that $\cF_{T''}(R)$ is not
well-defined for $D(n+1|n)$. For instance, for $n>1$,  for each $k>0$
one has $v(-2k\vareps_1)=0$  so
$\max\supp\bigl(t_{-2k\vareps_1}(R)\bigr)=0$ and the sum
$\sum_{k=1}^{\infty} t_{-2k\vareps_1}(R)$ is not well-defined; hence
$\cF_{T''}(R)$ is not well-defined as well.
\end{rem}

\subsection{Step 2}\label{pfsuppU}
By Section~\ref{infprod}, $\hat{R}$ is an invertible element of $\cR'$. 
From representation theory we know that since $\fhg$ admits a 
Casimir element~\cite{Kbook}, Chapter II, the character of the trivial
$\fhg$-module is a linear combination of the characters of 
Verma $\fhg$-modules $M(\lambda)$, where $\lambda\in-\hat{Q}$ are such that
$(\lambda+\hat{\rho},\lambda+\hat{\rho})=(\hat{\rho},\hat{\rho})$.
Since the character of $M(\lambda)$ is equal to $\hat{R}^{-1}e^{\lambda}$,
we obtain
$$1=\sum_{\substack
{\lambda\in \hat{Q}^-,\\ 
(\lambda+\hat{\rho},\lambda+\hat{\rho})
=(\hat{\rho},\hat{\rho})}} a_{\lambda}\hat{R}^{-1}e^{\lambda},$$
where $a_{\lambda}\in\mathbb{Z}$. This can be rewritten as
$$\hat{R}e^{\hat{\rho}}=\sum_{\substack{\lambda\in \hat{\rho}-\hat{Q}^+,\\ 
(\lambda,\lambda)
=(\hat{\rho},\hat{\rho})}} a_{\lambda}e^{\lambda},$$
that is $\supp (\hat{R})\subset U$, see~(\ref{defU}) for notation.

It remains to verify the inclusion 
$\supp\bigl(\cF_{T'}(Re^{\hat{\rho}})\bigr)\subset U$.
The denominator identity for $\fg$ (see~\cite{KW},\cite{G1}) takes the form
$$Re^{{\rho}}=
\cF_{W''}\bigl(\frac{e^{\rho}}{\prod_{\beta\in S}(1+e^{-\beta})}\bigr),$$
where $S:=\{\vareps_i-\delta_i\}_{i=1}^n$ (the identity for $\fgl(n|n)$
immediately follows from the  identity for $\fsl(n|n)$).
Since $\rho=\hat{\rho}$ is $\hat{W}$-invariant, this implies
$$t_{\mu}(Re^{\hat{\rho}})=e^{\hat{\rho}}\sum_{w\in W''}\sgn(w)
\prod_{\beta\in S}(1+e^{-t_{\mu}w\beta})^{-1}.$$
For each $t_{\mu}\in T'$ and $w\in W''$ one has 
$$\supp\bigl(\prod_{\beta\in S}(1+e^{-t_{\mu}w\beta})^{-1}  \bigr)\subset V,
\text{ where }
V:=\mathbb{Z}\{ t_{\mu}w\beta|\ \beta\in S\}\cap \hat{Q}^-.$$
Since $(t_{\mu}w\beta,t_{\mu}w\beta')=(\beta,\beta')=
(t_{\mu}w\beta,\hat{\rho})=(\hat{\rho},\beta)=0$
for any $\beta,\beta'\in S$ , one has 
$(V,V)=(V,\hat{\rho})=0$. Therefore $V+\hat{\rho}\subset U$
so $\supp\bigl( t_{\mu}(Re^{\hat{\rho}})\bigr)\subset U$ for each $\mu$.
This establishes the required inclusion
$\supp\bigl(\cF_{T'}(Re^{\hat{\rho}})\bigr)\subset U$
and completes the proof of~(\ref{suppUU}).

\subsection{Step 3}\label{Y1Y2}
Let us  deduce the inclusion 
$\supp(\hat{R}^{-1}e^{\hat{\rho}}\cdot{\cF_{T'}(Re^{\hat{\rho}})} )
\subset (\hat{Q}^-)^{\hat{W}}$ from~(\ref{suppUU}).

\subsubsection{}
\begin{lem}{lemai}
For any  simple finite-dimensional contragredient Lie superalgebra $\fg$
the term $\cF_{T'}(Re^{\hat{\rho}})$ is a $\hat{W}'$-anti-invariant element 
of $\cR_{\hat{W}'}$.
\end{lem}
\begin{proof}
In the light of Section~\ref{compex}, it is enough to present 
$\cF_{T'}(Re^{\hat{\rho}})$ in the form $\cF_{\hat{W}'}(Y)$ for
some $Y\in\cR_{\hat{W}}$.  Let $R_0', R_0''$ be the 
Weyl denominators for $\Delta',\Delta''$ respectively
(i.e., $R_0'=\prod_{\alpha\in\Delta'_+}(1-e^{-\alpha})$). Below
we will prove the formula
\begin{equation}\label{R0''}
\cF_{T'}(Re^{\hat{\rho}})=
 \cF_{\hat{W}'}\bigl(\frac{R_0''e^{\hat{\rho}}}{R_1}\bigr).
\end{equation}
By Section~\ref{infprod}, $R_1^{-1}R_0''e^{\hat{\rho}}\in\cR_{\hat{W}}$,
so the formula establishes the required assertion.

Let us show that the right-hand side of~(\ref{R0''}) is  well-defined.
Since $R_0''$ is $\hat{W}'$-invariant, it is enough to verify that
$\cF_{\hat{W}'}\bigl(e^{\hat{\rho}}R_1^{-1}\bigr)$ is a well-defined element 
of $\cR$. For $\fg\not=\fgl(n|n), D(n+1|n)$ 
this is proven in~\cite{G}, 2.4.1 (i).
Consider the case  $\fg=\fgl(n|n), D(n+1|n)$. Since $\hat{\rho}$
is $\hat{W}$-invariant, it is enough to check that
$\cF_{\hat{W}'}(R_1^{-1})$ is a well-defined element 
of $\cR$. By Section~\ref{cR'}, for each $w\in \hat{W}'$  one has
$$\max\supp\bigl(w(R_1^{-1})\bigr)=
\sum_{\beta\in\Delta_{1+}: w\beta<0}w\beta.$$
In particular, $\supp\bigl(w(R_1^{-1})\bigr)\subset\hat{Q}^-$, so,
if  the sum $\cF_{\hat{W}'}(R_1^{-1})
=\sum_{w\in\hat{W}'}\sgn w\cdot w(R_1^{-1})$ is well-defined,
it lies in $\cR$. In order to see that this sum
is  well-defined let us check that for each $\nu\in\hat{Q}^-$ the set
$$X(\nu):=\{w\in \hat{W}'|\ \sum_{\beta\in\Delta_{1+}: w\beta<0}
w\beta\geq \nu\}$$
is finite. One has
$$X(\nu)\subset \{w\in \hat{W}'|\ 
\forall \beta\in\Delta_{1+}\ w\beta\geq \nu\}.$$
Write $\nu=-k\delta+\nu'$, where 
$k\geq 0,\ \nu'\in Q$, and write $w\in X(\nu)$ in the the form
 $w=t_{\mu}y$, where $t_{\mu}\in T',
y\in W'$. Since $w\beta=y\beta-(y\beta,\mu)\delta$ for $\beta\in\Delta_{1+}$,
one has $(y\beta,\mu)\geq -k$ for each $\beta\in\Delta_{1+}$.
Since
$\{\vareps_i-\delta_i,\delta_i-\varesp_{i+1}\}\subset\Delta_{1+}$, this gives
 $|(\mu,y\delta_i)|\leq k$  for $i=1,\ldots,n$.
Combining the facts that $W'$ is a subgroup of signed permutation
of $\{\delta_j\}_{j=1}^n$ and that $(\mu,\delta_i)$ is integral for each $i$,
we conclude that $X(\nu)$ is finite. 
Thus $\cF_{\hat{W}'}\bigl(\frac{R_0''}{R_1}\bigr)$ is a well-defined
element of $\cR$.

Now let us prove the formula~(\ref{R0''}).
Recall that $\rho=\rho'_0+\rho''_0-\rho_1$, where
$$\rho'_0:=\sum_{\alpha\in\Delta'_{0+}}\alpha/2,\ \ 
\rho''_0:=\sum_{\alpha\in\Delta''_{0+}}\alpha/2,\ \ 
\rho_1:=\sum_{\beta\in\Delta_{1+}}\beta/2.$$

The Weyl denominator identity for $\Delta''_0$ takes the form
$$R_0'e^{\rho'_0}=\cF_{W'}(e^{\rho'_0}).$$
Since
$R_1e^{\rho_1}=\prod_{\beta\in\Delta_{1+}}(e^{\beta/2}+e^{-\beta/2})$ 
is $W$-invariant and $R_0''e^{\rho''_0}$ is $W'$-invariant, we get
$$Re^{{\rho}}=\frac{R_0''e^{\rho_0''}}{R_1e^{\rho_1}}
\cdot\cF_{W'}(e^{\rho'_0})=
\cF_{W'}\bigl(\frac{e^{\rho'_0}R_0''e^{\rho_0''}}{R_1e^{\rho_1}}\bigr)=
\cF_{W'}\bigl(\frac{R_0''e^{\rho}}{R_1}\bigr).$$
Using the $W$-invaraince of $\hat{\rho}-\rho$, we obtain
$$\cF_{T'}\bigl(Re^{\hat{\rho}})=
\cF_{T'}\bigl(\cF_{W'}\bigl(\frac{R_0''e^{\hat{\rho}}}{R_1}\bigr)\bigr)=
\cF_{\hat{W}'}\bigl(\frac{R_0''e^{\hat{\rho}}}{R_1}\bigr)$$
as required. This completes the proof.
\end{proof}

\subsubsection{}
\begin{prop}{ai}
Let $\fg$ be a simple
finite-dimensional contragredient Lie superalgebra. One has
$$\supp(\hat{R}^{-1}e^{\hat{\rho}}\cdot{\cF_{T'}(Re^{\hat{\rho}})})
\subset (\hat{Q}^-)^{\hat{W}}=\hat{Q}^-\cap\hat{Q}^{\perp}.$$
\end{prop}
\begin{proof}
By Section~\ref{welldef}, $\cF_{T'}(Re^{\hat{\rho}})\in\cR$; 
by Section~\ref{infprod}, $\hat{R}^{-1}\in\cR$
so
$$Y:=\hat{R}^{-1}e^{-\hat{\rho}}\cdot \cF_{T'}(Re^{\hat{\rho}})\in\cR.$$

The affine root system $\hat{\Delta}'$ is a subsystem of $\hat{\Delta}_0$.
Set $\hat{\Delta}'_+=\hat{\Delta}'\cap \hat{\Delta}_+$ and
let $\hat{\Pi}'$ be the corresponding set of simple roots.
Fix $\hat{\rho}'\in\fhh^*$ such that $2(\hat{\rho}',\alpha)=(\alpha,\alpha)$
for each $\alpha\in\hat{\Pi}'$.

It is easy to see that $\hat{R}_0e^{\hat{\rho}'},
\hat{R}e^{\hat{\rho}}$ are $\hat{W}'$-anti-invariant elements
of $\cR'$ (see, for instance,~\cite{G}, 1.5.1).
Thus $\hat{R}_1e^{\hat{\rho}'-\hat{\rho}}=\hat{R}_0e^{\hat{\rho}'}\cdot
(\hat{R}e^{\hat{\rho}})^{-1}$ is a $\hat{W}'$-invariant
element of $\cR'$.
By Section~\ref{infprod}, $\hat{R}_1\in\cR_{\hat{W}}$
so  $\hat{R}_1e^{\hat{\rho}'-\hat{\rho}}$ is a $\hat{W}'$-invariant
element of $\cR_{\hat{W}}$. Using~\Lem{lemai}, we get
\begin{equation}\label{R0Y}
\hat{R}_0e^{\hat{\rho}'}Y=\hat{R}_1e^{\hat{\rho}'-\hat{\rho}}\cF_{T'}(R)
\ \text{ 
 is a $\hat{W}'$-anti-invariant element of }\cR_{\hat{W}'}.
\end{equation}

Write $Y=Y_1+Y_2$, where
$\supp(Y_1)=\supp(Y)\cap \hat{Q}^{\perp}$ and
$\supp(Y_2)=\supp(Y)\setminus \hat{Q}^{\perp}$.
Note that $Y_1,Y_2\in\cR$. Assume that $Y_2\not=0$. 
Let $\mu$ be a maximal element in $\supp (Y_2)$.
One has $\supp(\hat{R}^{-1})\subset \hat{Q}^-$ and
$\supp\bigl(\cF_{T'}(R)e^{\hat{\rho}}\bigr)
\subset \hat{\rho}-\hat{Q}^+$, by Section~\ref{cR'} 
and~(\ref{suppUU}) respectively.
Thus $\supp(Y)\subset \hat{Q}^-$ and so $\mu\in \hat{Q}^-$.

Since $\supp(Y_1)\subset\hat{Q}^{\perp}$,
$Y_1$ is a $\hat{W}$-invariant element of $\cR_{\hat{W}}$
so $\hat{R}_0e^{\hat{\rho}'}Y_1$
is a $\hat{W}'$-anti-invariant element of $\cR_{\hat{W}'}$. 
In the light of~(\ref{R0Y}), the product $\hat{R}_0e^{\hat{\rho}'}Y_2$
is also a $\hat{W}'$-anti-invariant element of $\cR_{\hat{W}'}$. 
Clearly, $\hat{\rho}'+\mu$ is a maximal element in
the support of $\hat{R}_0e^{\hat{\rho}'}Y_2$.
By Section~\ref{compex}, this support is the union of $\hat{W}'$-regular
orbits (recall that regularity means
that each element has the trivial stabilizer in $\hat{W}'$),
so  $\hat{\rho}'+\mu$ is a maximal element
in a regular $\hat{W}'$-orbit and thus
$\frac{2(\hat{\rho}'+\mu,\alpha)}{(\alpha,\alpha)}\not\in\mathbb{Z}_{\leq 0}$
for each  $\alpha\in\hat{\Pi}'$. 
Since $\mu\in\hat{Q}^-$ one has
$\frac{2(\mu,\alpha)}{(\alpha,\alpha)}\in\mathbb{Z}$ for each 
$\alpha\in\hat{\Pi}'$. Taking into account that 
$\frac{2(\hat{\rho}',\alpha)}{(\alpha,\alpha)}=1$ for each 
$\alpha\in\hat{\Pi}'$, we obtain
\begin{equation}\label{mualpha}
\forall\alpha\in\hat{\Pi}'\ \ 
\frac{2(\mu,\alpha)}{(\alpha,\alpha)}\in\mathbb{Z}_{\geq 0}.
\end{equation}
Recall
that $\delta=\sum_{\alpha\in\hat{\Pi}'}k_{\alpha}\alpha$ for some
$k_{\alpha}\in\mathbb{Z}_{>0}$ (see~\cite{Kbook}, Chapter VI). Since 
$\mu\in\hat{Q}^-$ one has $(\mu,\delta)=0$. Combining with~(\ref{mualpha}),
we get $(\mu,\alpha)=0$ for each $\alpha\in\hat{\Pi}'$ so
$\mu\in(\hat{\Delta}')^{\perp}$.

One has 
$$(\hat{\Delta}')^{\perp}\cap\hat{Q}=(\hat{Q}^{\perp}\cap\hat{Q})
\oplus V,$$
where the restriction of $(-,-)$ to $\mathbb{Q}V$ is negatively definite;
more precisely, one has
$$\begin{tabular}{|l || l | l| l|l|}
\hline
$\fg$ & $\fgl(n|n)$ &  $\fgl(m|n),\ m\not=n$ & 
\ \ \ $C(n)$ &  other cases\\
\hline\hline
$\hat{Q}^{\perp}\cap\hat{Q}$ & 
 $\mathbb{Z}\{\delta,\str\}$ & $\mathbb{Z}\delta$ &
$\mathbb{Z}\delta$ &$\mathbb{Z}\delta$\\
\hline
$V$ & $\mathbb{Z}\Delta''$ & $\mathbb{Z}\Delta''\oplus \mathbb{C}\xi$
&  $\mathbb{Z}\Delta''\oplus \mathbb{C}\xi$ & $\mathbb{Z}\Delta''$\\
\hline
\end{tabular}$$
For $\fg=\fgl(m|n),\ m\not=n$ and $\fg=C(n)$ the element $\xi$ is
given in~\cite{G}, 3.2;  one has
$(\Delta'',\xi)=0,\ (\xi,\xi)<0$. Since $V\subset \hat{Q}$,
one has $(V,\hat{Q}^{\perp})=0$.
Now combining the formulas
$\mu\in(\hat{Q}^{\perp}\cap\hat{Q})\oplus V,\ (\mu,\mu)=0$
with the fact that $(\nu,\nu)<0$ for each non-zero $\nu\in V$,
we obtain $\mu\in\hat{Q}^{\perp}\cap\hat{Q}=\hat{Q}^{\hat{W}}$, 
which contradicts to the construction
of $Y_2$. Hence $Y_2=0$ as required.
\end{proof}

\subsubsection{}\label{Corfq}
Using the  table in the proof of~\Prop{ai}, we obtain the following
corollary.

\begin{cor}{corfq}
For $\fg\not=\fgl(n|n)$ one has 
$f(q)\cdot\hat{R}e^{\hat{\rho}}=\cF_{T'}(Re^{\hat{\rho}})$
for some $f(q)=\sum_{k=0}^{\infty} a_kq^k$ ($a_k\in\mathbb{Z}$).
For $\fg=\fgl(n|n)$ one has 
$f(q,e^{\str})\cdot\hat{R}e^{\hat{\rho}}=\cF_{T'}(Re^{\hat{\rho}})$
for some $f(q,e^{\str})=\sum_{k=0}^{\infty}\sum_{m=-\infty}^{\infty} 
a_{k,m}q^ke^{m\cdot\str}$ ($a_{k,m}\in\mathbb{Z}$).
\end{cor}

\subsection{Step 4 for $\fg=\fgl(n|n),D(n+1|n)$}\label{ev}
In this subsection we complete the proof of the 
denominator identities~(\ref{denom0}) by proving the formulas~(\ref{fqintr}).
We prove them by taking a suitable evaluation of 
$\hat{R}^{-1}\sum_{t\in T'} t(R)$. 
By~\Cor{corfq}, 
$\hat{R}^{-1}\sum_{t\in T'} t(R)$ is equal to $f(q)$
for $D(n+1|n)$ and to $f(q,e^{\str})$ for $\fgl(n|n)$.
Now we consider $q$ as a real parameter between $0$ and $1$.
We choose the evaluation in such a way that
the evaluation of $\hat{R}^{-1}\sum_{t\in T'} t(R)$
is equal to the evaluation of $\hat{R}^{-1}R$.
 As a result, $f(q)$ (resp.,  $f(q,e^{\str})$) 
is equal to  the evaluation of $\hat{R}^{-1}R$,
which can be easily computed.

\subsubsection{Case $D(n+1|n)$}
Take a complex parameter $x$ and consider the following evaluation:
$e^{-\vareps_i}:=x^{a_i},\ e^{-\delta_j}:=-x^{b_j}$,
where $a_i,\ (i=1,\ldots,n+1),\ b_j, (j=1,\ldots,n)$ are integers such that 
$a_i\pm b_j\not=0, a_i\pm a_j\not=0,
b_i\pm b_j\not=0, b_i\not=0$ for all indexes $i,j$. We denote
the evaluation of $R$ (resp., $\hat{R}$) by $R(x)$
(resp., $\hat{R}(x)$). The functions $R(x),\hat{R}(x)$ are meromorphic.
One has
$$R(x)=\frac{\prod_{1\leq i<j\leq n+1}(1-x^{a_i\pm a_j})
\cdot\prod_{1\leq i<j\leq n}(1-x^{b_i\pm b_j})
\cdot\prod_{1\leq i\leq n}(1-x^{2b_i})}
{\prod_{1\leq i\leq j\leq n}(1-x^{a_i\pm b_j})
\prod_{1\leq j<i\leq n+1}(1-x^{a_i\pm b_j})}.$$
One readily sees that  $R(x)$ has a pole at $x=1$ of order $|\Delta_{1+}|-
|\Delta_{0+}|=n$.

One has
$$\left.\frac{\hat{R}(x)}{R(x)}\right|_{x=1}=
\frac{((1-q)^{\infty}_{q})^{\dim\fg_0}}
{((1-q)^{\infty}_{q})^{\dim\fg_1}}=((1-q)^{\infty}_{q})^
{\dim\fg_0-\dim\fg_1}=(1-q)^{\infty}_{q}.$$
In particular, $\hat{R}(x)$ also has a pole of order $n$  at $x=1$.

The evaluation of $(t_{\sum k_i\delta_i}(R))(x)$ is
$$
\frac{\prod_{1\leq i<j\leq n+1}(1-x^{a_i\pm a_j})\cdot\prod_{1\leq i\leq n}
(1-q^{-2k_i}x^{2b_i})\cdot \prod_{1\leq i<j\leq n}
(1-q^{-k_i\mp k_j}x^{b_i\pm b_j})}
{\prod_{1\leq i\leq j\leq n}(1-q^{\mp k_j}x^{a_i\pm b_j})
\prod_{1\leq j<i\leq n+1}(1-q^{\mp k_j}x^{-a_i\pm b_j})}$$
which is a meromorphic function.
Let $s$ be the number of zeros among $k_1,\ldots,k_n$. Then at $x=1$ 
the order of zero of the numerator is at least
is $n(n+1)+s^2$, and the order of zero of the denominator is
$2(n+1)s$. Therefore  at $x=1$ the function
$(t_{\sum k_i\delta_i}(R))(x)$ has the pole of order at most
$2(n+1)s-n(n+1)-s^2=n+1-(n+1-s)^2$; in particular,
$(t_{\sum k_i\delta_i}(R))(x)$ has the pole of order 
at most $n$ and it is equal to $n$ iff $n=s$ that is $\sum k_i\delta_i=0$
and $(t_{\sum k_i\delta_i}(R))(x)=R(x)$.

We conclude that 
$(\hat{R}(x))^{-1}\cdot\sum_{t\in T': t\not=\id} (t(R))(x)$ is  holomorphic
 at $x=1$ and its value is equal to zero, and that
$(\hat{R}(x))^{-1}\cdot\sum_{t\in T'} (t(R))(x)$ is  holomorphic
 at $x=1$ and its value is equal to $\frac{R(x)}{\hat{R}(x)}|_{x=1}$.
In the light of~\Cor{corfq} we obtain  
$$f(q)=\left.\frac{R(x)}{\hat{R}(x)}\right|_{x=1}=((1-q)^{\infty}_{q})^{-1}.$$

\subsubsection{Case $\fgl(n|n)$}
Fix $y>1$.
Take a complex parameter $x$ and consider the following evaluation
$$e^{-\vareps_1}:=y,\ e^{-\vareps_i}:=x^{i},\text{ for } i=2,\ldots,n
\ e^{-\delta_i}:=-x^{-i}\text{ for } i=1,\ldots,n.$$
The functions $R(x), \hat{R}(x)$  are meromorphic.
One has
$$R(x)=\frac{\prod_{1<i\leq n}(1-yx^{-i})\cdot
\prod_{1<i<j\leq n}(1-x^{i-j})\cdot\prod_{1\leq i<j\leq n}(1-x^{j-i})}
{\prod_{1\leq i\leq n}(1-yx^{i})\cdot\prod_{1<i\leq j\leq n}(1-x^{i+j})
\cdot\prod_{1\leq j<i\leq n}(1-x^{-i-j})}.$$
Therefore the function $R(x)$ has a pole of order $n-1$ at $x=1$.

One has
$$\left.\frac{\hat{R}(x)}{R(x)}\right|_{x=1}=
\frac{((1-q)^{\infty}_{q})^{\dim\fg_0-2(n-1)}
\cdot ((1-qy)^{\infty}_{q})^{n-1}\cdot ((1-qy^{-1})^{\infty}_{q})^{n-1}}
{((1-q)^{\infty}_{q})^{\dim\fg_1-2n}\cdot 
((1-qy)^{\infty}_{q})^{n}\cdot ((1-qy^{-1})^{\infty}_{q})^{n}}.$$
Thus $\hat{R}(x)$ also has a pole of order $n-1$ at $x=1$.
Since $\dim\fg_0=\dim\fg_1$ and $e^{\str}=(-1)^ny^{-1}$ 
for $x=1$ we obtain
$$\left.\frac{\hat{R}(x)}{R(x)}\right|_{x=1}=\frac{((1-q)^{\infty}_{q})^2}
{(1-q(-1)^ne^{\str})^{\infty}_{q}\cdot (1-q(-1)^ne^{-\str})^{\infty}_{q}}.$$

One has
$$(t_{\sum k_i\delta_i}(R))(x,y)=\frac{\prod_{1<i\leq n}(1-yx^{-i})\cdot
\prod_{1<i<j\leq n}(1-x^{i-j})
\cdot\prod_{1\leq i<j\leq n}(1-q^{k_j-k_i}x^{j-i})}
{\prod_{1\leq i\leq n}(1-q^{k_i}yx^{i})\cdot\prod_{1<i\leq j\leq n}(1-
q^{k_j}x^{i+j})\cdot\prod_{1\leq j<i\leq n}(1-q^{-k_j}x^{-i-j})},$$
which is a meromorphic function. 

Let $s$ be the number of zeros among $k_1,\ldots,k_n$. Then at $x=1$ 
the order of zero of the numerator is at least 
$\frac{(n-1)(n-2)+s(s-1)}{2}$, and the order of zero of the denominator is 
$(n-1)s$. Therefore  at $x=1$ the function
$(t_{\sum k_i\delta_i}(R))(x,y)$ has the pole of order at most
$(n-1)s-\frac{(n-1)(n-2)+s(s-1)}{2}=\frac{3n-s-2-(n-s)^2}{2}$, so
the order is at most $n-1$ and it is equal to $n-1$ iff $s=n-1,n$.
Notice that $s\not=n-1$, since $\sum k_i=0$. Therefore
the pole has order $n-1$ iff $\sum k_i\delta_i=0$.

We conclude that the function $(\hat{R}(x))^{-1}(\cF_{T'}(R))(x)$ 
is holomorphic  at $x=1$ and its value is equal to 
$\frac{R(x)}{\hat{R}(x)}|_{x=1}$. Using~\Cor{Corfq} we obtain  
$$f(q,e^{\str})=\left.\frac{R(x)}{\hat{R}(x)}\right|_{x=1}=
\frac{(1-q(-1)^ne^{\str})^{\infty}_{q}
\cdot(1-q(-1)^ne^{-\str})^{\infty}_{q}}{((1-q)^{\infty}_{q})^2}.$$

\subsection{Step 4 for $\fg\not=\fgl(n|n),D(n+1|n), D(2,1,a)$}
\label{other}
In this case the dual Coxeter number is non-zero. Recall that $q=e^{-\delta}$.
Write $f(q)=\sum_{k=0}^{\infty} a_ke^{-k\delta}$.
Since $f(q)\cdot\hat{R}e^{\hat{\rho}}=\cF_{T'}(Re^{\hat{\rho}})$, we have
$$\sum_{k=1}^{\infty} a_ke^{-\delta}\cdot \hat{R}e^{\hat{\rho}}=
\cF_{T'}(Re^{\hat{\rho}})-a_0\hat{R}e^{\hat{\rho}}.$$
By~(\ref{suppUU}), for any $\nu$ in the support of
the right-hand side, one has $(\nu,\nu)=(\hat{\rho},\hat{\rho})$,
and  for any $\nu$ in the support of
the left-hand side one has $(\nu,\nu)=(\hat{\rho},\hat{\rho})-
2k(\delta,\hat{\rho})$ for some $k>0$. Since $(\hat{\rho},\delta)$
is equal to the dual Coxeter number, which is non-zero, we conclude
that the intersection of supports is empty. Hence $f(q)=a_0$. Since
the coefficient of $e^{\hat{\rho}}$ in $\hat{R}e^{\hat{\rho}}$ is equal to one,
$a_0$ is equal to the coefficient of $e^{\hat{\rho}}$ in 
$\cF_{T'}(Re^{\hat{\rho}})$. As it is shown in~\cite{G}, this coefficient
is equal to one so $f(q)=1$ as required.

\section{Other forms of denominator identity}
Recall that denominator identity for a basic Lie superalgebra can be written 
in the form
\begin{equation}\label{findenom}
Re^{\rho}=\cF_{W^{\sharp}}\bigl(\frac{e^{\rho}}
{\prod_{\beta\in S}(1+e^{-\beta})}\bigr),
\end{equation}
where $W^{\sharp}:=W'$ for $\fg\not=D(n+1|n), D(2,1,a)$ and
$W^{\sharp}:=W''$ for $\fg=D(n+1|n), D(2,1,a)$, and
$S\subset \Pi$ is the maximal isotropic system (see~\cite{KW},\cite{G1}).
If the dual Coxeter number of $\fg$ is non-zero
the affine denominator identity for $\fg$ can be written 
in the form
$$
\hat{R}e^{\hat{\rho}}=\cF_{\hat{W^{\sharp}}}\bigl(\frac{e^{\hat{\rho}}}
{\prod_{\beta\in S}(1+e^{-\beta})}\bigr)
$$
see~\cite{KW},\cite{G}. In this section we will show
that for $\fgl(n|n)$ the denominator identity can be written 
in a similar form:
\begin{equation}\label{slnnS}
\hat{R}e^{\rho}=f(q,e^{\str})\cdot\cF_{\hat{W'}}\bigl(\frac{e^{\rho}}
{\prod_{\beta\in S}(1+e^{-\beta})}\bigr),
\end{equation}
and that the  denominator identities
for $D(n+1|n)$  can not be written in a similar form,
since the expressions $\cF_{\hat{W''}}\bigl(\frac{e^{\rho}}
{\prod_{\beta\in S}(1+e^{-\beta})}\bigr)$,
$\cF_{\hat{W'}}\bigl(\frac{e^{\rho}}
{\prod_{\beta\in S}(1+e^{-\beta})}\bigr)$ are not well defined.

\subsection{Case $D(n+1|n)$}\label{DD}
Let us show that the expressions $\cF_{\hat{W''}}\bigl(\frac{e^{\rho}}
{\prod_{\beta\in S}(1+e^{-\beta})}\bigr)$,
$\cF_{\hat{W'}}\bigl(\frac{e^{\rho}}
{\prod_{\beta\in S}(1+e^{-\beta})}\bigr)$ are not well-defined
for $D(n+1|n)$. Fix $\Pi$ as in Section~\ref{rootsys} and
recall that $\rho=0$.

We repeat the reasonings of Section~\ref{welldef}. One has
$$\sum_{\beta\in V(w)}w\beta\in
\supp\bigl(\frac{1}
{\prod_{\beta\in S}(1+e^{-w\beta})}\bigr)
\subset \sum_{\beta\in V_S(w)}w\beta-\hat{Q}^+
\subset\hat{Q}^-, $$
where
$$V_S(w)=\{\beta\in S| w\beta<0\}.$$

Therefore $1\in \supp\bigl(\frac{1}
{\prod_{\beta\in S}(1+e^{-w\beta})}\bigr)$ iff $wS\subset\Delta_+$.

Take $S=\{\vareps_i-\delta_i\}$; then $t_{\mu}S\subset\Delta_+$
if $(\vareps_i-\delta_i,\mu)<0$ for all $i$ which holds
for all  $\mu\in\sum\mathbb{Z}_{<0}\vareps_i$ and all
$\mu\in\sum\mathbb{Z}_{>0}\delta_i$. Hence 
the sums $\cF_{\hat{W''}}\bigl(\frac{e^{\rho}}
{\prod_{\beta\in S}(1+e^{-\beta})}\bigr)$,
$\cF_{\hat{W'}}\bigl(\frac{e^{\rho}}
{\prod_{\beta\in S}(1+e^{-\beta})}\bigr)$ contain infinitely many summands
equal to $1$ and thus they are not  well-defined.

\subsection{Case $\fgl(n|n)$}
Fix $\Pi$ as in Section~\ref{rootsys}; then $S=\{\vareps_i-\delta_i\}$.

In order to deduce the formula~(\ref{slnnS}) from~(\ref{findenom}) and
(\ref{denom0}) it is enough to verify that the expression
$$\cF_{\hat{W'}}\bigl(\frac{e^{\rho}}{\prod_{\beta\in S}(1+e^{-\beta})}\bigr)
= e^{\rho}\cF_{\hat{W'}}
\bigl(\frac{1}{\prod_{\beta\in S}(1+e^{-\beta})}\bigr)$$
is well-defined (since $\rho=\str/2$ is $\hat{W}$-invariant).
As in Section~\ref{welldef}, it amounts to show that
$$X_S(\nu):=\{w\in\hat{W}'|\ \sum_{\beta\in V_S(w)} w\beta\geq -\nu\}$$
is finite for any $\nu\in\hat{Q}^+$ 
(where $V_S(w)$ is defined as in Section~\ref{DD}).
As in Section~\ref{welldef}, writing $\nu=k\delta+\nu_+$, where 
$\nu_+\in\mathbb{Z}\Delta$, we get
$$X_S(\nu)\subset 
\{t_{\mu}y|\mu\in T', y\in W'\text{ s.t. } (yS,\mu)\geq -k\}.$$
Since $y$ permutes $\delta_i$s, $t_{\mu}y\in X_S(\nu)$ forces
$(\delta_i,\mu)\geq -k$ for all $i$. Taking into account
that $\mu$ lies in the $\mathbb{Z}$-span of $\delta_i$ and
$(\mu,\sum_{i=1}^n\delta_i)=0$, we conclude that $X_S(\nu)$ is finite.
This establishes~(\ref{slnnS}).


\end{document}